\documentclass[a4paper, 12pt]{article}

\usepackage{fullpage}
\usepackage{authblk}
\usepackage{amsmath}
\usepackage{amssymb}
\usepackage{gensymb}
\usepackage{xfrac}
\usepackage{url}
\usepackage{graphicx}
\usepackage{csquotes}

\MakeOuterQuote{"}

\title{Kullback-Leibler Divergence \\ for the Normal-Gamma Distribution}
\author[1,3]{Joram Soch}
\author[1,2]{Carsten Allefeld}
\affil[1]{Bernstein Center for Computational Neuroscience, Berlin, Germany}
\affil[2]{Berlin Center for Advanced Neuroimaging, Berlin, Germany}
\affil[3]{Department of Psychology, Humboldt-Universität zu Berlin, Germany}
\date{}

\begin{document}

\setcounter{page}{0}
\vspace*{1em}
\begin{center}
	\LARGE
		Kullback-Leibler Divergence \\ for the Normal-Gamma Distribution \\ \vspace{1em}
	\large
		Joram Soch\textsuperscript{1,3,\textbullet} \& Carsten Allefeld\textsuperscript{1,2}
\end{center}
\begin{flushleft}
	\normalsize
		\hspace{2.5em}
		\textsuperscript{1} Bernstein Center for Computational Neuroscience, Berlin, Germany \\
		\hspace{2.5em}
		\textsuperscript{2} Berlin Center for Advanced Neuroimaging, Berlin, Germany \\
		\hspace{2.5em}
		\textsuperscript{3} Department of Psychology, Humboldt-Universit\"{a}t zu Berlin, Germany \\ \vspace{1em}
		\hspace{2.5em}
		\textbullet \ Corresponding author: \url{joram.soch@bccn-berlin.de}.		
\end{flushleft}
\vspace*{1em}

\begin{abstract}
\noindent
We derive the Kullback-Leibler divergence for the normal-gamma distribution and show that it is identical to the Bayesian complexity penalty for the univariate general linear model with conjugate priors. Based on this finding, we provide two applications of the KL divergence, one in simulated and one in empirical data.
\end{abstract}

\vspace*{1em}
\tableofcontents

\pagebreak
\section{Introduction}

Let $x$ be a $k \times 1$ random vector and $y > 0$ be a random variable. Then, $x$ and $y$ are said to follow a \textit{normal-gamma distribution} (NG distribution), if their joint probability density function is given by

\begin{equation} \label{eq:NG-pdf1}
p(x,y) = \mathrm{N}(x; \mu, (y \Lambda)^{-1}) \cdot \mathrm{Gam}(y; a, b)
\end{equation}

where $\mathrm{N}(x; \mu, \Sigma)$ denotes a multivariate normal density with mean $\mu$ and covariance $\Sigma$ and $\mathrm{Gam}(x; a, b)$ denotes a gamma density with shape $a$ and rate $b$. In full, the density function is given by (Koch, 2007, p.~55)

\begin{equation} \label{eq:NG-pdf2}
p(x,y) = \sqrt{\frac{|y \Lambda|}{(2 \pi)^k}} \, \exp\left[ -\frac{y}{2} (x-\mu)^T \Lambda (x-\mu) \right] \cdot \frac{{b}^{a}}{\Gamma(a)} \, y^{a-1} \, \exp[-b y] \; .
\end{equation}

\vspace{0.5em}

The \textit{Kullback-Leibler divergence} (KL divergence) is a non-symmetric distance measure for two probability distributions $P$ and $Q$ and is defined as

\begin{equation} \label{eq:disc-KL}
\mathrm{KL}[P||Q] = \sum_{i \in \Omega} P(i) \, \ln \frac{P(i)}{Q(i)} \; .
\end{equation}

For continuous probability distribtions $P$ and $Q$ with probability density functions $p(x)$ and $q(x)$ on the same domain $X$, it is given by (Bishop, 2006, p.~55)

\begin{equation} \label{eq:cont-KL}
\mathrm{KL}[P||Q] = \int_{X} p(x) \, \ln \frac{p(x)}{q(x)} \, \mathrm{d}x \; .
\end{equation}

The KL divergence becomes important in information theory and statistical inference. Here, we derive the KL divergence for two NG distributions with vector-valued $x$ and real-positive $y$ and provide two examples of its application.

\pagebreak
\section{Theory}

\subsection{Multivariate normal KL divergence}

First, consider two multivariate normal distributions over the $k \times 1$ vector $x$ specified by

\begin{equation}
\begin{split}
p(x) &= \mathrm{N}(x; \mu_1, \Sigma_1) \\
q(x) &= \mathrm{N}(x; \mu_2, \Sigma_2)
\end{split}
\end{equation}

According to equation (\ref{eq:cont-KL}), the KL divergence of $P$ from $Q$ is defined as

\begin{equation}
\mathrm{KL}[P||Q] = \int_{\mathbb{R}^k} \mathrm{N}(x; \mu_1, \Sigma_1) \, \ln \frac{\mathrm{N}(x; \mu_1, \Sigma_1)}{\mathrm{N}(x; \mu_2, \Sigma_2)} \, \mathrm{d}x \; .
\end{equation}

Using the multivariate normal density function

\begin{equation} \label{eq:mvn-pdf}
\mathrm{N}(x; \mu, \Sigma) = \frac{1}{\sqrt{(2 \pi)^n |\Sigma|}} \, \exp\left[ -\frac{1}{2} (x-\mu)^T\Sigma^{-1}(x-\mu) \right] \; ,
\end{equation}

it evaluates to (Duchi, 2014)

\begin{equation} \label{eq:mvn-KL}
\mathrm{KL}[P||Q] = \frac{1}{2} \left[ (\mu_2 - \mu_1)^T \Sigma_2^{-1} (\mu_2 - \mu_1) + \mathrm{tr}(\Sigma_2^{-1} \Sigma_1) - \ln \frac{|\Sigma_1|}{|\Sigma_2|} - k \right] \; .
\end{equation}

\subsection{Univariate gamma KL divergence}

Next, consider two univariate gamma distributions over the real-positive $y$ specified by

\begin{equation}
\begin{split}
p(y) &= \mathrm{Gam}(y; a_1, b_1) \\
q(y) &= \mathrm{Gam}(y; a_2, b_2)
\end{split}
\end{equation}

According to equation (\ref{eq:cont-KL}), the KL divergence of $P$ from $Q$ is defined as

\begin{equation}
\mathrm{KL}[P||Q] = \int_{0}^{\infty} \mathrm{Gam}(y; a_1, b_1) \, \ln \frac{\mathrm{Gam}(y; a_1, b_1)}{\mathrm{Gam}(y; a_2, b_2)} \, \mathrm{d}y \; .
\end{equation}

Using the univariate gamma density function

\begin{equation} \label{eq:gam-pdf}
\mathrm{Gam}(y; a, b) = \frac{{b}^{a}}{\Gamma(a)} \, y^{a-1} \, \exp[-b y] \quad \text{for} \quad y > 0 \; ,
\end{equation}

it evaluates to (Penny, 2001)

\begin{equation} \label{eq:gam-KL}
\mathrm{KL}[P||Q] = a_2 \, \ln \frac{b_1}{b_2} - \ln \frac{\Gamma(a_1)}{\Gamma(a_2)} + (a_1 - a_2) \, \psi(a_1) - (b_1 - b_2) \, \frac{a_1}{b_1}
\end{equation}

where $\psi(x)$ is the digamma function.

\pagebreak
\subsection{Normal-gamma KL divergence}

Now, consider two normal-gamma distributions over $x$ and $y$ specified by

\begin{equation}
\begin{split}
p(x,y) &= \mathrm{N}(x; \mu_1, (y \Lambda_1)^{-1}) \, \mathrm{Gam}(y; a_1, b_1) \\
q(x,y) &= \mathrm{N}(x; \mu_2, (y \Lambda_2)^{-1}) \, \mathrm{Gam}(y; a_2, b_2)
\end{split}
\end{equation}

According to equation (\ref{eq:cont-KL}), the KL divergence of $P$ from $Q$ is defined as

\begin{equation} \label{eq:NG-KL0}
\mathrm{KL}[P||Q] = \int_{0}^{\infty} \int_{\mathbb{R}^k} p(x,y) \, \ln \frac{p(x,y)}{q(x,y)} \, \mathrm{d}x \, \mathrm{d}y \; .
\end{equation}

Using the law of conditional probability, it can be evaluated as follows:

\begin{equation} \label{eq:NG-KL1}
\begin{split}
\mathrm{KL}[P||Q] &= \int_{0}^{\infty} \int_{\mathbb{R}^k} p(x|y) p(y) \, \ln \frac{p(x|y) p(y)}{q(x|y) q(y)} \, \mathrm{d}x \, \mathrm{d}y \\
&= \int_{0}^{\infty} p(y) \int_{\mathbb{R}^k} p(x|y) \, \ln \frac{p(x|y)}{q(x|y)} \, \mathrm{d}x \, \mathrm{d}y \\
&+ \int_{0}^{\infty} p(y) \, \ln \frac{p(y)}{q(y)} \int_{\mathbb{R}^k} p(x|y) \, \mathrm{d}x \, \mathrm{d}y \\
&= \left\langle \mathrm{KL}[p(x|y)||q(x|y)] \right\rangle_{p(y)} + \mathrm{KL}[p(y)||q(y)]
\end{split}
\end{equation}

In other words, the KL divergence for two normal-gamma distributions over $x$ and $y$ is equal to the sum of a multivariate normal KL divergence regarding $x$ conditional on $y$, expected over $y$, and a univariate gamma KL divergence regarding $y$. Together with equation (\ref{eq:mvn-KL}), the first term becomes

\begin{equation}
\begin{split}
&\left\langle \mathrm{KL}[p(x|y)||q(x|y)] \right\rangle_{p(y)} \\
&= \left\langle \frac{1}{2} \left[ (\mu_2 - \mu_1)^T (y \Lambda_2) (\mu_2 - \mu_1) + \mathrm{tr}\left( (y \Lambda_2) (y \Lambda_1)^{-1} \right) - \ln \frac{|(y \Lambda_1)^{-1}|}{|(y \Lambda_2)^{-1}|} - k \right] \right\rangle_{p(y)} \\
&= \left\langle \frac{y}{2} (\mu_2 - \mu_1)^T \Lambda_2 (\mu_2 - \mu_1) + \frac{1}{2} \, \mathrm{tr}(\Lambda_2 \Lambda_1^{-1}) - \frac{1}{2} \ln \frac{|\Lambda_2|}{|\Lambda_1|} - \frac{k}{2} \right\rangle_{p(y)} \; .
\end{split}
\end{equation}

Using the relation $y \sim \mathrm{Gam}(a,b) \Rightarrow \left\langle y \right\rangle = a/b$, we have

\begin{equation} \label{eq:exp-mvn-KL}
\begin{split}
\left\langle \mathrm{KL}[p(x|y)||q(x|y)] \right\rangle_{p(y)} = \frac{1}{2} \frac{a_1}{b_1} (\mu_2 - \mu_1)^T \Lambda_2 (\mu_2 - \mu_1) + \frac{1}{2} \, \mathrm{tr}(\Lambda_2 \Lambda_1^{-1}) - \frac{1}{2} \ln \frac{|\Lambda_2|}{|\Lambda_1|} - \frac{k}{2} \; .
\end{split}
\end{equation}

Thus, from (\ref{eq:exp-mvn-KL}) and (\ref{eq:gam-KL}), the KL divergence in (\ref{eq:NG-KL1}) becomes

\begin{equation} \label{eq:NG-KL2}
\begin{split}
\mathrm{KL}[P||Q] &= \frac{1}{2} \frac{a_1}{b_1} \left[ (\mu_2 - \mu_1)^T \Lambda_2 (\mu_2 - \mu_1) \right] + \frac{1}{2} \, \mathrm{tr}(\Lambda_2 \Lambda_1^{-1}) - \frac{1}{2} \ln \frac{|\Lambda_2|}{|\Lambda_1|} - \frac{k}{2} \\
&+ a_2 \, \ln \frac{b_1}{b_2} - \ln \frac{\Gamma(a_1)}{\Gamma(a_2)} + (a_1 - a_2) \, \psi(a_1) - (b_1 - b_2) \, \frac{a_1}{b_1} \; .
\end{split}
\end{equation}

\pagebreak
\subsection{The Bayesian model evidence}

Consider Bayesian inference on data $y$ using model $m$ with parameters $\theta$. In this case, Bayes' theorem is a statement about the posterior density:

\begin{equation} \label{eq:BT}
p(\theta|y,m) = \frac{p(y|\theta,m) \, p(\theta|m)}{p(y|m)} \; .
\end{equation}

The denominator $p(y|m)$ acts as a normalization constant on the posterior density $p(\theta|y,m)$ and according to the law of marginal probability is given by

\begin{equation} \label{eq:ME1}
p(y|m) = \int p(y|\theta,m) \, p(\theta|m) \, \mathrm{d}\theta \; .
\end{equation}

This is the probability of the data given only the model, regardless of any particular parameter values. It is also called "marginal likelihood" or "model evidence" and can act as a model quality criterion in Bayesian inference, because parameters are integrated out of the likelihood.

For computational reasons, only the logarithmized or log model evidence (LME) $\mathrm{L}(m) = \ln p(y|m)$ is of interest in most cases. By rearranging equation (\ref{eq:BT}), the model evidence can be represented as

\begin{equation} \label{eq:ME2}
p(y|m) = \frac{p(y|\theta,m) \, p(\theta|m)}{p(\theta|y,m)} \; .
\end{equation}

Logarithmizing both sides of the equation and taking the expectation with respect to the posterior density over model parameters $\theta$ gives the LME

\begin{equation} \label{eq:LME1}
\mathrm{L}(m) = \int p(\theta|y,m) \, \ln p(y|\theta,m) \, \mathrm{d}\theta - \int p(\theta|y,m) \, \ln \frac{p(\theta|y,m)}{p(\theta|m)} \, \mathrm{d}\theta \; .
\end{equation}

Using this reformulation, the LME as a model quality measure can be naturally decomposed into an accuracy term, the posterior expected likelihood, and a complexity term, the KL divergence between the posterior and the prior distribution:

\begin{equation} \label{eq:LME2}
\begin{split}
\mathrm{L}(m) &= \mathrm{Acc}(m) - \mathrm{Com}(m) \\
\mathrm{Acc}(m) &= \left\langle \log p(y|\theta,m) \right\rangle_{p(\theta|y,m)} \\
\mathrm{Com}(m) &= \mathrm{KL} \left[ p(\theta|y,m) || p(\theta|m) \right]
\end{split}
\end{equation}

Intuitively, the accuracy acts increasing and the complexity acts decreasing on the log model evidence. This reflects the capability of the LME to select models that achieve the best balance between accuracy and complexity, i.e.~models that explain the observations sufficiently well (high accuracy) without employing too many principles (low complexity). The fact that the complexity term is a KL divergence between posterior and prior means that models with prior assumptions that are close to the posterior evidence receive a low complexity penalty, because one is not surprised very much when accepting such a model which renders the Bayesian complexity a measure of surprise.

\pagebreak
\subsection{The general linear model}

Consider multiple linear regression using the univariate general linear model (GLM)

\begin{equation} \label{eq:GLM}
y = X \beta + \varepsilon, \; \varepsilon \sim N(0,\sigma^2 V)
\end{equation}

where $y$ is an $n \times 1$ vector of measured data, $X$ is an $n \times p$ matrix called the design matrix, $\beta$ is a $p \times 1$ vector of weight parameters called regression coefficients and $\varepsilon$ is an $n \times 1$ vector of errors or noise. These residuals are assumed to follow a multivariate normal distribution whose covariance matrix is the product of a variance factor $\sigma^2$ and an $n \times n$ correlation matrix $V$. Usually, $X$ and $V$ are known while $\beta$ and $\tau$ are unknown parameters to be inferred via model estimation.

For mathematical convenience, we rewrite $\sigma^2 = 1/\tau$ and $V = P^{-1}$ so that equation (\ref{eq:GLM}) implies the following likelihood function:

\begin{equation} \label{eq:GLM-LF}
p(y|\beta,\tau) = \mathrm{N}(y; X \beta, (\tau P)^{-1}) \; .
\end{equation}

The conjugate prior relative to this likelihood function is a normal-gamma distribution on the model parameters $\beta$ and $\tau$ (Koch, 2007, ch. 2.6.3):

\begin{equation} \label{eq:GLM-NG-prior}
\begin{split}
p(\beta|\tau) &= \mathrm{N}(\beta; \mu_0, (\tau \Lambda_0)^{-1}) \\
p(\tau) &= \mathrm{Gam}(\tau; a_0, b_0)
\end{split}
\end{equation}

Due to the conjugacy of (\ref{eq:GLM-NG-prior}) to (\ref{eq:GLM-LF}), the posterior is also a normal-gamma distribution

\begin{equation} \label{eq:GLM-NG-post}
\begin{split}
p(\beta|\tau,y) &= \mathrm{N}(\beta; \mu_n, (\tau \Lambda_n)^{-1}) \\
p(\tau|y) &= \mathrm{Gam}(\tau; a_n, b_n)
\end{split}
\end{equation}

where the posterior parameters in (\ref{eq:GLM-NG-post}) are given by (Koch, 2007, ch. 4.3.2)

\begin{equation} \label{eq:GLM-NG-post-par}
\begin{split}
\mu_n &= \Lambda_n^{-1} (X^T P y + \Lambda_0 \mu_0) \\
\Lambda_n &= X^T P X + \Lambda_0 \\
a_n &= a_0 + \frac{n}{2} \\
b_n &= b_0 + \frac{1}{2} (y^T P y + \mu_0^T \Lambda_0 \mu_0 - \mu_n^T \Lambda_n \mu_n)
\end{split}
\end{equation}

From (\ref{eq:LME2}), the complexity for the model defined by (\ref{eq:GLM-LF}) and (\ref{eq:GLM-NG-prior}) is given by

\begin{equation} \label{eq:GLM-NG-Com1}
\mathrm{Com}(m) = \mathrm{KL} \left[ p(\beta,\tau|y) || p(\beta,\tau) \right] \; .
\end{equation}

In other words, the complexity penalty for a general linear model with normal-gamma priors (GLM-NG) is identical to a KL divergence between two NG distributions and using (\ref{eq:NG-KL2}) can be written in terms of the prior and posterior parameters as

\begin{equation} \label{eq:GLM-NG-Com2}
\begin{split}
\mathrm{Com}(m) &= \frac{1}{2} \frac{a_n}{b_n} \left[ (\mu_n - \mu_0)^T \Lambda_0 (\mu_n - \mu_0) \right] + \frac{1}{2} \, \mathrm{tr}(\Lambda_0 \Lambda_n^{-1}) - \frac{1}{2} \ln \frac{|\Lambda_0|}{|\Lambda_n|} - \frac{p}{2} \\
&+ a_0 \, \ln \frac{b_n}{b_0} - \ln \frac{\Gamma(a_n)}{\Gamma(a_0)} + (a_n - a_0) \, \psi(a_n) - (b_n - b_0) \, \frac{a_n}{b_n} \; .
\end{split}
\end{equation}

\pagebreak
\section{Application}

\subsection{Polynomial basis functions}

Consider a linear model with polynomial basis functions (Bishop, 2006, p.~5) given by

\begin{equation} \label{eq:PBF1}
y = \sum_{i=0}^p c_i \, x^i + \varepsilon \quad \text{with} \quad -1 \leq x \leq +1 \; .
\end{equation}

Essentially, this model assumes that $y$ is an additive mixture of polynomial terms $x^i$ weighted with the coefficients $c_i$ with $i = 1,\ldots,p$ where the natural number $p$ is called the model order. This means that $p = 0$ corresponds to a constant value (plus noise $\varepsilon$); $p = 1$ corresponds to a linear function; $p = 2$ corresponds to a quadratic pattern; $p = 3$ corresponds to a 3rd degree polynomial etc.

Given that $x$ is an $n \times 1$ vector of real numbers between $-1$ and $+1$, this model can be rewritten as a GLM given in equation (\ref{eq:GLM}) with

\begin{equation} \label{eq:PBF2}
X = 
\begin{bmatrix}
x_1^0  & x_1^1  & \cdots & x_1^p  \\
x_2^0  & x_2^1  & \cdots & x_2^p  \\
\vdots & \vdots & \ddots & \vdots  \\
x_n^0  & x_n^1  & \cdots & x_n^p
\end{bmatrix}
\quad \text{and} \quad
\beta = 
\begin{bmatrix}
c_0 \\ c_1 \\ \vdots \\ c_p
\end{bmatrix} \; .
\end{equation}

Based on this reformulation, we simulate polynomial data. We perform $N = 100$ simulations with $n = 100$ data points in each simulation. We generate simulated data based on a true model order $p_{\mathrm{true}} = 5$ and analyze these data using a set of models ranging from $p_{\mathrm{min}} = 0$ to $p_{\mathrm{max}} = 20$. The predictor $x$ is equally spaced between $-1$ and $+1$ and design matrices $X_p$ are created according to equation (\ref{eq:PBF2}).

In each simulation, six regression coefficients $\beta_{\mathrm{true}}$ are drawn independently from the standard normal distribution $\mathrm{N}(0,1)$. Then, Gaussian noise $\varepsilon$ is sampled from the multivariate normal distribution $\mathrm{N}(0,\sigma_{\varepsilon}^2 I_n)$ with a residual variance of $\sigma_{\varepsilon}^2 = 1$. Finally, simulated data is generated as as $y = X_5 \beta_{\mathrm{true}} + \varepsilon$.

Then, for each $p \in \left\lbrace 0,\ldots,20 \right\rbrace$, Bayesian model estimation is performed using the design matrix $X_p$, a correlation matrix $V = I_n$ and the prior distributions (\ref{eq:GLM-NG-prior}) with the prior parameters $\mu_0 = 0_p$, $\Lambda_0 = I_p$ invoking a standard multivariate normal distribution and $a_0 = 1$, $b_0 = 1$ invoking a relatively flat gamma prior. Posterior parameters are calculated using equation (\ref{eq:GLM-NG-post-par}) and give rise to the model complexity via (\ref{eq:GLM-NG-Com2}) as well as model accuracy and the log model evidence via (\ref{eq:LME2}).

Average LMEs, accuracies and complexities are shown in Figure~1. One can see that the true model order is correctly identified by the maximal log model evidence. This is achieved by an increasing complexity penalty which outweighs the saturating accuracy gain for models with $p > 5$. This demonstrates that the KL divergence for the NG distribution can be used to select polynomial basis functions when basis sets cannot be separated based on model accuracy alone.

\pagebreak
\begin{center} \label{fig:Figure_PBF}
\includegraphics[width=0.99\linewidth, clip=true, trim=0 100 0 75]{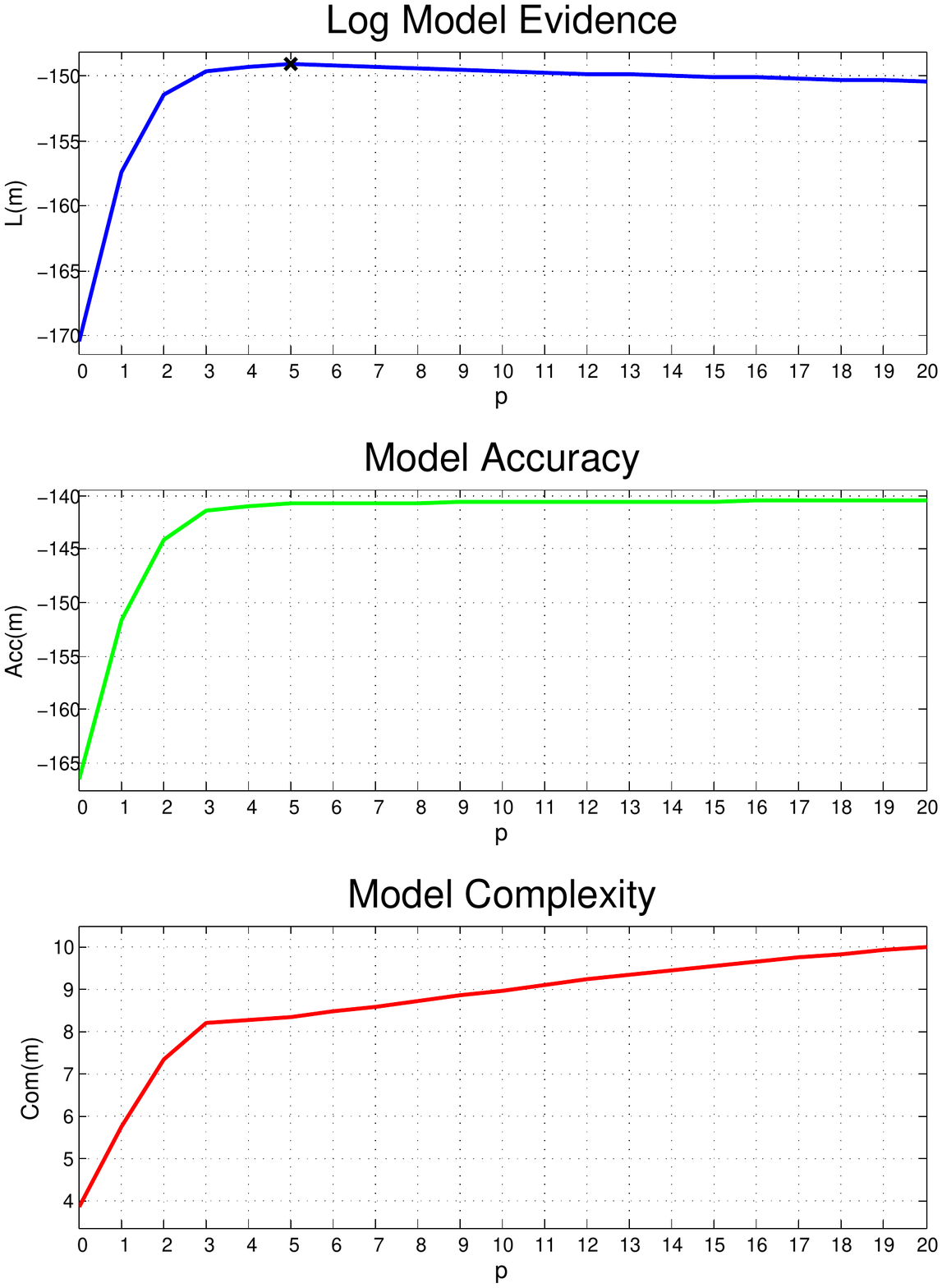}
\end{center}

\textbf{Figure 1.} Bayesian model selection for polynomial basis functions. All displays have model order on the x-axis and average model quality measures (across simulations) on the y-axis. Intuitively, the model accuracy (middle panel) increases with model order, but saturates at around $-140$ with no major increase after $p = 5$. Moreover, the model complexity (lower panel) -- which is the KL divergence between prior and posterior distribution -- also grows with model order, but switches to a linear increase at around $p = 5$ reaching a value of $10$ at $p = 20$. Together, this has the consequence that the log model evidence (upper panel) is maximal for $p = 5$ (black cross) where exact values are: $\mathrm{Acc}(m) = -140.77, \, \mathrm{Com}(m) = 8.34, \, \mathrm{L}(m) = -149.11$.

\pagebreak
\subsection{Neuroimaging model selection}

In neuroimaging, especially functional magnetic resonance imaging (fMRI), GLMs as given by equation (\ref{eq:GLM}) are applied to time series of neural data $y$ (Friston et al., 1995). The design matrix $X$ is specified by the temporal occurrence of experimental conditions and the covariance structure $V$ is estimated from residual auto-correlations. Model estimation and statistical inference are performed "voxel-wise", i.e.~separately for each measurement location in the brain, usually referred to as the "mass-univariate GLM".

Here, we analyze data from a study on orientation pop-out processing (Bogler et al., 2013). During the experimental paradigm, the screen showed a $3 \times 7$ array of homogeneous bars oriented either 0\degree{}, 45\degree{}, 90\degree{} or 135\degree{} relative to the vertical axis. This background stimulation changed every second and was interrupted by trials in which one target bar on the left and one target bar on the right were independently rotated either 0\degree{}, 30\degree{}, 60\degree{} or 90\degree{} relative to the rest of the stimulus display. Those trials of orientation contrast (OC) lasted 4 seconds and were alternated with inter-trial intervals of 7, 10 or 13 seconds. Each combination of OC on the left side and OC on the right side was presented three times resulting in 48 trials in each of the 5 sessions lasting 672 seconds.

After fMRI data preprocessing (slice-timing, realignment, normalization, smoothing), two different models of hemodynamic activation were applied to the fMRI data. The first model (GLM I) considers the experiment a factorial design with two factors (left vs. right OC) having four levels (0\degree{}, 30\degree{}, 60\degree{}, 90\degree{}). This results in $4 \times 4 = 16$ possible combinations or experimental conditions modelled by $16$ onset regressors convolved with the canonical hemodynamic response function (HRF). The second model (GLM II) puts all trials from all conditions into one HRF-convolved regressor and encodes orientation contrast using a parametric modulator (PM) that is given as $\mathrm{PM} = \mathrm{deg}/{90\degree{}}$ with $\mathrm{deg} = ( 0\degree{}, 30\degree{}, 60\degree{}, 90\degree{} )$, resulting in $\mathrm{PM} = ( 0, \sfrac{1}{3}, \sfrac{2}{3}, 1 )$, such that the parametric modulator is proportional to orientation contrast. There was one PM for each factor of the design, i.e.~one PM for left OC and one PM for right OC. Note that both models encode the same information and that every signal that can be identified using GLM II can also be characterized using GLM I, but not vice versa, because the first model allows for a greater flexibility of activation patterns across experimental conditions than the second.

For these two models, we performed Bayesian model estimation. To overcome the challenge of having to specify prior distributions on the model parameters, we apply cross-validation across fMRI sessions. This gives rise to a cross-validated log model evidence (cvLME) as well as cross-validated accuracies and complexities for each model in each subject. We then performed a paired t-test to find voxels where GLM II has a significantly higher cvLME than GLM I. Due to the specific assumptions in GLM II and the higher flexibility of GLM I, we assumed that these differences might be primarily based on a complexity advantage of GLM II over GLM I.

We focus on visual area 4 (V4) that is known to be sensitive to orientation contrast. Within left V4, specified by a mask from a separate localizer paradigm (Bogler et al., 2013), we identified the peak voxel ($[\text{x y z}] = [-15, -73, -5]$~mm) defined by the maximal t-value ($t = \mathrm{3.95}$) and extracted log model evidence as well as model accuracy and model complexity from this voxel for each subject. Differences in LME, accuracy and complexity are shown in Figure~2. Again, the model complexity enables a model selection that would not be possible based on the model accuracy alone.

\pagebreak
\begin{center} \label{fig:Figure_NMS}
\includegraphics[width=0.99\linewidth, clip=true, trim=0 100 0 75]{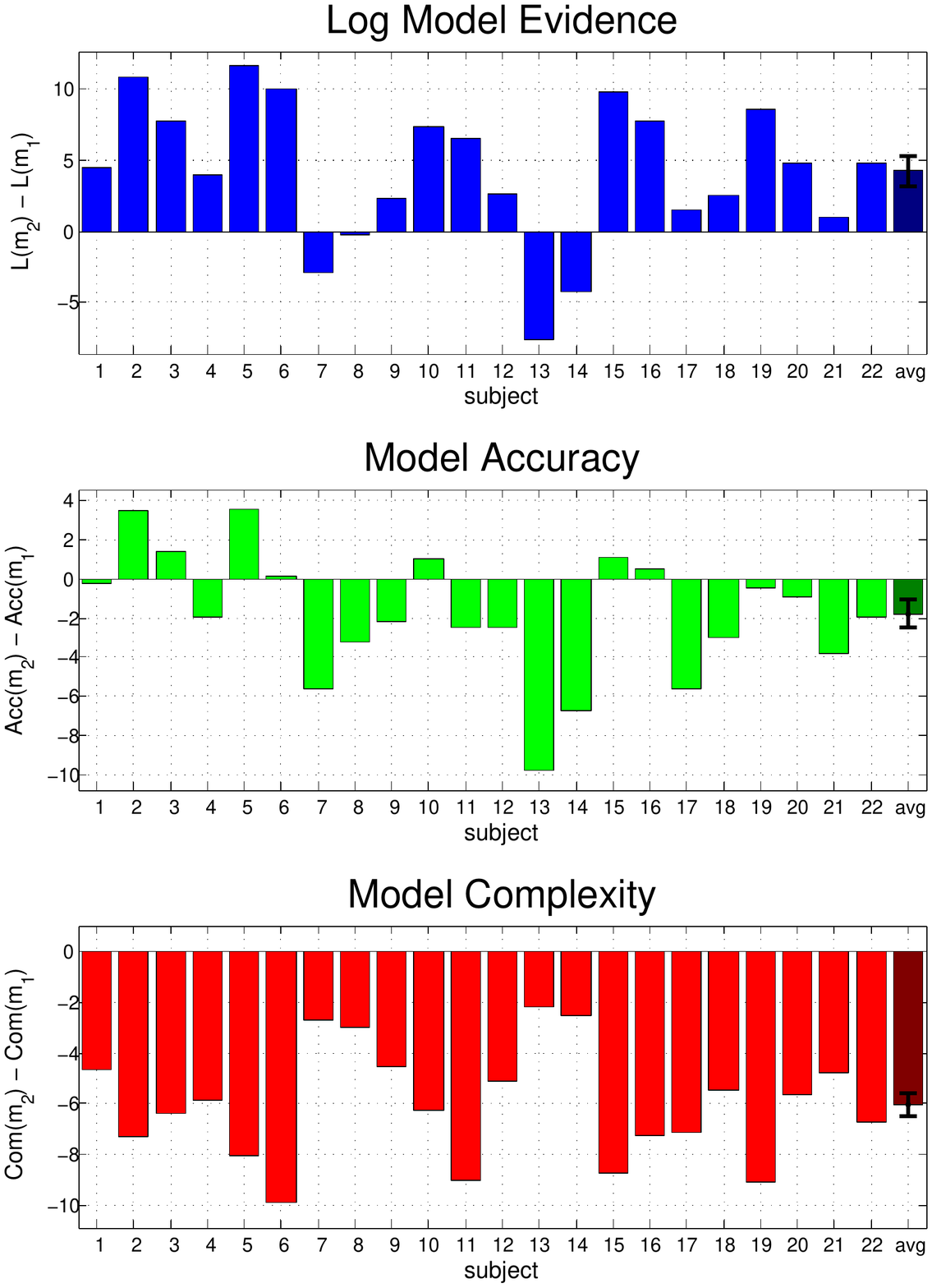}
\end{center}

\textbf{Figure 2.} Bayesian model selection for orientation pop-out processing. All displays have subject on the x-axis and difference in model qualities ($\mathrm{m}_1$ = GLM I, $\mathrm{m}_2$ = GLM II) on the y-axis. Interestingly, there is a slight disadvantage for GLM II regarding only the model accuracy (middle panel), its mean difference across subjects being smaller than zero. However, model complexity (lower panel) -- measured as the KL divergence between prior and posterior distribution -- is consistently higher for GLM I. Together, this has the consequence that the log model evidence (upper panel) most often favors GLM II. Average values are: $\Delta \mathrm{Acc} = -1.78, \, \Delta \mathrm{Com} = -6.02, \, \Delta \mathrm{L} = 4.24$.

\pagebreak
\section{Conclusion}

We have derived the Kullback-Leibler divergence of two normal-gamma distributions using earlier results on the KL divergence for multivariate normal and univariate gamma distributions. Moreover, we have shown that the KL divergence for the NG distribution occurs as the complexity term in the univariate general linear model when using conjugate priors.

Analysis of simulated and empirical data demonstrates that the complexity penalty has the desired theoretical features, namely to quantify the relative informational content of two generative models and to detect model differences that cannot be detected by just relying on model accuracy, e.g.~given by the maximum log-likelihood (as in information criteria like AIC or BIC) or the posterior log-likelihood (as in the Bayesian log model evidence).

\section{References}

\renewcommand{\section}[2]{}

\end{document}